\documentclass[runningheads,a4paper]{llncs}

\usepackage{amssymb}
\usepackage{graphicx}

\usepackage{amsmath}
\usepackage{algpseudocode}
\usepackage{algorithm}

\usepackage{url}
\urldef{\mailsa}\path|{d.andreagiovanni,enrico.natalizio}@hds.utc.fr|
\urldef{\mailsb}\path|nardin.an@gmail.com|
\newcommand{\keywords}[1]{\par\addvspace\baselineskip
\noindent\keywordname\enspace\ignorespaces#1}

\begin{document}

\mainmatter

\title{A fast ILP-based Heuristic for the robust design of Body Wireless Sensor Networks\thanks{
This is the authors' final version of the paper published in G. Squillero and K. Sim (Eds.): EvoApplications 2017, Part I, LNCS 10199, pp. 1-17, 2017.
DOI: 10.1007/978-3-319-55849-3\_16
The final publication is available at Springer via http://dx.doi.org/10.1007/978-3-319-55849-3\_16}
}

\author{
Fabio D'Andreagiovanni$^{1}$
\and
Antonella Nardin$^{2}$
\and
Enrico Natalizio$^{1}$
}

\pagestyle{empty}


\institute{
$^{1}$Sorbonne Universit\'es, Universit\'e de Technologie de Compi\`egne, CNRS, Heudiasyc UMR 7253, CS 60319, 60203 Compi\`egne, France\\
$^{2}$Universit\`a degli Studi Roma Tre,
    Via Ostiense 169, 00154 Roma, Italy\\
\mailsa\\
\mailsb\\
}

\toctitle{Lecture Notes in Computer Science}
\tocauthor{Authors' Instructions}
\maketitle

\begin{abstract}
We consider the problem of optimally designing a body wireless sensor network, while taking into account the uncertainty of data generation of biosensors. Since the related min-max robustness Integer Linear Programming (ILP) problem can be difficult to solve even for state-of-the-art commercial optimization solvers, we propose an original heuristic for its solution. The heuristic combines deterministic and probabilistic variable fixing strategies, guided by the information coming from strengthened linear relaxations of the ILP robust model, and includes a very large neighborhood search for reparation and improvement of generated solutions, formulated as an ILP problem solved exactly. Computational tests on realistic instances show that our heuristic finds solutions of much higher quality than a state-of-the-art solver and than an effective benchmark heuristic.
\end{abstract}

\section{Introduction}

\noindent
A \emph{Wireless Sensor Network} (WSN) can be essentially described as a network of typically small and portable wireless devices, the sensors, which are spread on an area to collect data in a cooperative way and then forward the data to one or more collectors, commonly called \emph{sinks}.
Recently, the application of WSNs in healthcare has received a lot of attention and, just to cite two major examples, WSNs have been used to monitor the health conditions of patients in hospitals and to remotely monitor people under health risk when they are at home \cite{KoEtAl07}.

In this work, we focus attention on a topic related to healthcare applications of WSNs: the design of \emph{body area networks}.
A \emph{Body Area Network} (BAN) is a WSN where wireless sensors (\emph{biosensors}) are placed over or inside the body of a person to collect biomedical data. The biosensors generate data and transmit them to one or more sinks for storing or processing.  For a detailed introduction to BANs, we refer the reader to the works \cite{ChEtAl10,MiEtAl12}.
Designing a BAN essentially consists in deciding the topology of the network and how data are routed from the biosensors to the sinks. This constitutes a classical WSN design problem (see e.g., \cite{GuEtAl11,NaLoVi08}). However, since BANs are deployed on human bodies, their design need particular attention and present specific challenges that are not shared with other WSNs design problems \cite{ChEtAl10,NeJeBe16}.  A critical question is in particular constituted by the peculiar high-loss propagation behaviour of wireless signals through and over the human body: in contrast to canonical wireless networks, where high losses can be handled by increasing power emissions (see e.g., \cite{DA15,DAMaSa13,GeEtAl16}), in BANs power emissions must be contained to both avoid damages to human tissues, due to overheating, and to preserve the charge of sensor batteries, whose substitution can result very uncomfortable for patients. Controlling energy consumption is thus a major aim in BAN design and is typically achieved through multi-hop routing, implemented through relay nodes, which are wireless devices that act as intermediate nodes between sinks and sensors and allow transmission of reduced power over shorter distances \cite{El14,TsEtAl12,YoEtAl15},

Nowadays there is a rich literature about BANs, in particular about technical aspects concerning the definition of energy-efficient routing protocols and the study of the peculiar propagation condition in human bodies \cite{ChEtAl10,NeJeBe16}. In contrast, there is still a limited amount of work devoted to the design of BANs in terms of optimization models and algorithms. This fact has been highlighted also in the two relevant previous works \cite{El14,DANa15}: in \cite{El14}, the design problem of a BAN is formulated as a mixed integer linear program where multi-path routing and relay deployment is established in order to minimize the total cost of deployment of a BAN;  \cite{DANa15} instead investigates a robust optimization model for tackling the data generation uncertainty and proposes a Mixed Integer Programming (MIP)-based heuristic for the solution of the resulting challenging optimization problem.

In this work, we consider a \emph{scenario-based min-max} robust optimization model for the design of BANs that takes into account the uncertainty of data generation of BAN sensors.
Our main original contribution is a new ILP heuristic for solving the robust design problem, based on combining deterministic and probabilistic variable fixing strategies
guided by peculiar linear relaxations of the robust optimization model. In comparison to \cite{DANa15}, our new algorithm does
not just fix the variables expressing routing decisions, but also employs an initial deterministic fixing phase of the variables modelling the activation of relay nodes.
Computational tests on realistic BAN instances show that our new heuristic produces solutions that are not only deeply better than those produced by a state-of-the-art optimization solver, but are also significantly better than those found by the algorithm of \cite{DANa15}.

\section{The Body Area Network Design Problem}
\label{sec:BAND}

In this section, we identify the elements of a BAN that are relevant for modelling purposes and we derive a network optimization model for the energy-efficient design of a BAN. For a more detailed description of mathematical optimization modelling for BANs, we refer the reader to \cite{DANa15} and \cite{El14}.
\\
\textbf{System Elements of a BAN.}
Any BAN can be described as a \emph{set of biosensors}, denoted by $B$, that produce biomedical data while monitoring a human body. The data are collected by a \emph{set of sinks}, denoted by $S$. The biosensors and the sinks are located in positions over or inside the human body that are usually precisely pre-established (for example, one can think about the electrodes of a Holter monitor, which must be positioned on specific spots of the chest of a patient for monitoring heart activity). For each sink $s \in S$, each biosensor $b \in B$ generates a volume of data $d_{bs} \geq 0$ (typically, a \emph{bitrate} measured in bit/s).

In order to improve energy efficiency, the biosensors do not transmit their data directly to the sink according to a single-hop direct communication, but rely on a multi-hop routing strategy. 
Multi-hop routing can be implemented through relay nodes, which have the task of receiving and forwarding the biomedical data. The positions of relay nodes is not fixed and can be chosen. Such positioning choice, if done wisely, can greatly improve the energy efficiency of the BAN. The optimization problem that we consider is indeed related to the optimal positioning of relays in order to minimize energy consumption. In such problem, we are given a set of candidate locations for the relays and an upper bound on the number of deployable relays and we must take two decisions: i) establish the number of deployed relays; ii) choose the location of the deployed relays.

We introduce a set $R$ to represent the potentially deployable relays. Each potential relay $r \in R$ is characterized by a unique position in/over the body and we must decide whether it is deployed or not. Each $r$ has also a capacity $c_r \geq 0$ that represents the maximum bitrate that it can manage.

The transmission of data from any BAN device (biosensor, relay node or sink) to another BAN device is based on a directional wireless link. As in \cite{DANa15,El14}, we assume that the devices employs a TDMA (Time Division Multiple Access) protocol, which allows the devices to transmit on the same frequency band without interfering. When either transmitting or receiving data, the BANs devices consume energy according to the following formulas \cite{BrEtAl07}:
\begin{eqnarray}
\label{energyFormula}
        \begin{array}{lll}
            E_{\text{TX}} (v,\delta) &=& E_{\text{TX$_{CIRC}$}} \cdot v
            \hspace{0.1cm} + \hspace{0.1cm}
            E_{\text{TX$_{AMP}$}}(\lambda) \cdot \delta^{\lambda} \cdot v
            \\
            E_{\text{RX}} (v) &=& E_{\text{RX$_{CIRC}$}} \cdot v
        \end{array}
\end{eqnarray}

\noindent
where $E_{\text{TX}}$ is the total transmission energy and $E_{\text{RX}}$ the total receiver energy (expressed in joules).
We remark that $E_{\text{TX}}, E_{\text{RX}}$ are a function of the volume of transmitted/received data $v$ (expressed in bits) and of the distance $\delta$ (expressed in meters) between the transmitter and the receiver. Additionally: $E_{\text{TX$_{CIRC}$}}$, $E_{\text{RX$_{CIRC}$}}$ are the energy consumed by the circuits to respectively transmit and receive a single bit; $E_{\text{TX$_{AMP}$}}(\lambda)$ is the energy consumed by the transmitting amplifier and $\lambda$ is the path loss exponent in the signal attenuation formula.
\\
\textbf{A flow-based Integer Linear Program for BAN design.}
It is natural to trace back the energy-efficient design of BAN to a network flow optimization problem. Specifically, we trace it back to a
variant of a Multicommodity Flow Problem (MCFP), a classical network flow problem, where the aim is to decide how to install routing capacity and how to route a set of commodities in a network, minimizing the total routing and installation cost, while not exceeding the capacity of installed network elements. For an introduction to capacitated network design and MCFPs, we refer the reader to \cite{Be98} and \cite{DAKrPu15}.

\noindent
The BAN can be naturally modelled through a directed graph $G(V,A)$ where:

\smallskip
\noindent
1) the set of \emph{vertices} $V$ contains one element for each wireless device of the network - biosensor, sink and relay. The set $V$ thus corresponds to the union of three disjoint sets of vertices:
i) the set $B$ of vertices corresponding to biosensors;
ii) the set $R$ of vertices corresponding to potentially deployable relays;
iii) the set $S$ of vertices corresponding to sinks.
        We therefore have $V = B \cup R \cup S$.

        Each BAN device can communicate with other devices that are within its transmission range. The transmission range may vary on the basis of the propagation conditions and of the transmission power of the device (see e.g., \cite{BrEtAl07}). We denote the subsets of devices that are within the transmission range of a device as follows:
a) for each biosensor $b \in B$, we distinguish the subsets $R_b \subseteq R$, and $S_b \subseteq S$ representing the relays and sinks within the range of $b$, respectively;\\
b) for each potential relay $r \in R$, we distinguish the subsets $R_r \subseteq R$, and $S_r \subseteq S$ representing the relays and sinks within the range of $r$, respectively;\\
c) more generally, given a vertex $i \in V$, representing any type of BAN device, we denote by $V_i \subseteq V$ the subset of vertices representing devices within the transmission range of $i$.\footnote{We note that we assume that each biosensor $b \in B$ never acts as a receiver and \emph{only generates and transmits} data. So we do not characterize the subsets $B_r, B_s \subseteq B$ of biosensors within the range of a relay $r$ or a sink $s$. Furthermore, we assume that each sink $s$ never acts as a transmitter and \emph{only receives} data. So we do not characterize the subsets $B_s \subseteq B$, $R_s \subseteq R$ of biosensors and relays within the range of a sink $s$´.}

\smallskip
\noindent
2) the set of \emph{arcs} $A$ contains one element for each wireless link that can be established between a pair of wireless devices. An arc $a = (i,j) \in A$ is an ordered pair of vertices that models a directional wireless link from a device $i \in V$ to another device $j \in V_i$ within the range of $i$. We respectively call tail and head the vertices $i,j$ of $a = (i,j)$. The set $A$ is the union of four disjoint sets of arcs:\\
a) the set $A_{B \rightarrow S}$ of arcs $(i,j)$ such that the tail is a biosensor and the head is a sink within the range of the biosensor, i.e. $i \in B$, $j \in S_i$. They represent transmissions of biomedical data directly to sinks;\\
b) the set $A_{B \rightarrow R}$ of arcs $(i,j)$ such that the tail is a biosensor and the head is a relay within the range of the biosensor, i.e. $i \in B$, $j \in R_i$. They represent transmissions from a biosensor to a relay;\\
c) the set $A_{R \leftrightarrow R}$ of arcs $(i,j)$ such that both the tail and the head are relay nodes, i.e. $i,j \in R$ with $j \in R_i$. They represent wireless links between relays;\\
d) the set $A_{R \rightarrow S}$ of arcs $(i,j)$ such that the tail is a relay node and the head is a sink within the range of the relay, i.e. $i \in R$, $j \in S_i$. They represent transmissions from a relay to a sink.
Therefore we have $A = A_{B \rightarrow S} \cup A_{B \rightarrow R} \cup A_{R \leftrightarrow R} \cup A_{R \rightarrow S}$.

We can now rewrite the energy formulas \eqref{energyFormula} in terms of the graph introduced above. When data are transmitted on an arc $a = (i,j) \in A$, the energy consumed is the sum of the energy consumed by $i$ to transmit and the energy consumed by $j$ to receive. The energy consumed to send one unit of data from $i$ to $j$ is:
\begin{eqnarray}
\label{energyCoefficient}
E_{ij}
\hspace{0.05cm}=\hspace{0.05cm}
            E_{\text{TX}} (1,\delta_{ij}) + E_{\text{RX}} (1)
\hspace{0.05cm}=\hspace{0.05cm}
            \left[E_{\text{TX$_{\text{CIRC}}$}}
            +
            E_{\text{TX$_{\text{AMP}}$}}(\lambda_{ij}) \cdot \delta_{ij}^{\lambda_{ij}}\right]
            + E_{\text{RX$_{\text{CIRC}}$}}
\end{eqnarray}

\noindent
which is obtained from formulas \eqref{energyFormula} for $v = 1$.
Here, $\lambda_{ij}$ and $\delta_{ij}$ are respectively the path loss coefficient and the distance between $i$ and $j$.

Using all the notation and elements introduced until now, we can formally state the BAN design problem that we consider, referring to the problem definition and modelling that has been initially provided in \cite{DANa15}:
\begin{definition}[The Body Area Network Design Problem - BAND]
Given: 1) a BAN modeled as a directed graph $G(V,A)$, where $V = B \cup R \cup S$ is the set of vertices and $A = A_{B \rightarrow S} \cup A_{B \rightarrow R} \cup A_{R \leftrightarrow R} \cup A_{R \rightarrow S}$ is the set of arcs; \;
2) the bitrate $d_{bs} \geq 0$ of data generated by each biosensor $b \in B$ for each sink $s \in S$;  \;
3) the capacity $c_r \geq 0$ of each relay $r \in R$; \;
4) the energy coefficients $E_{ij} \geq 0$ expressing the total energy consumed to send 1 data unit from $i$ to $j$; \;
\\
the BAND consists in choosing which relays are activated and which single-paths are used to route the flow of data generated by each biosensor for each sink, in order to minimize the total energy consumption.
\qed
\end{definition}

\noindent
The BAND is based on taking two decisions: 1) which relays to deploy and 2) which single-paths to use to route the data generated by each biosensor for each sink. These two decisions can be modeled through two families of binary decision variables:
1) binary \emph{relay deployment variables} $y_{r} \in \{0,1\} \hspace{0.2cm} \forall \hspace{0.1cm} r \in R$ such that $y_{r}$ equals 1 if relay $r$ is deployed and 0 otherwise;
2) binary \emph{unsplittable flow variables} $x^{bs}_{ij} \in \{0,1\} \hspace{0.2cm} \forall \hspace{0.1cm} b \in B, s \in S, (i,j) \in A$  such that $x^{bs}_{ij}$ equals 1
                                if all the data generated by biosensor $b$ for sink $s$ are routed on arc $(i,j)$ and 0 otherwise.

These variables are employed in the following ILP problem that we denote by BAND-ILP \cite{DANa15}:
\begin{align}
\min \hspace{0.1cm} & \sum_{b \in B} \sum_{s \in S} \sum_{(i,j) \in A}
\hspace{0.1cm} E_{ij} \hspace{0.1cm} d_{bs} \hspace{0.1cm} x^{bs}_{ij}
&&\mbox{(BAND-ILP)}
\nonumber
\\
&
- \sum_{\substack{(b,j) \in \\ A_{B \rightarrow R} \cup A_{B \rightarrow S}}} d_{bs} \; x^{bs}_{bj} \hspace{0.1cm} = - d_{bs} \;
&&
b \in B, s \in S
\label{BAND-BLP_bioConservation}
\\
&
\sum_{\substack{(j,r) \in \\ A_{B \rightarrow R} \cup A_{R \leftrightarrow R}}}
d_{bs} \; x^{bs}_{jr}
- \sum_{\substack{(r,j) \in \\ A_{R \leftrightarrow R} \cup A_{R \rightarrow S}}}
d_{bs} \; x^{bs}_{rj} = 0
&&
b \in B, s \in S, r \in R
\label{BAND-BLP_relayConservation}
\\
&
\sum_{\substack{(j,s) \in \\ A_{B \rightarrow S} \cup A_{R \rightarrow S}}} d_{bs} \; x^{bs}_{js} \hspace{0.1cm} = \hspace{0.1cm} d_{bs} \;
&&
b \in B, s \in S
\label{BAND-BLP_sinkConservation}
\\
&
\sum_{\substack{(r,j) \in \\ A_{R \leftrightarrow R} \cup A_{R \rightarrow S}}} d_{bs} \hspace{0.1cm} x^{bs}_{rj} \hspace{0.1cm} \leq \hspace{0.1cm} c_r \hspace{0.05cm} y_r
&&
r \in R
\label{BAND-BLP_capacity}
\\
&
\sum_{r \in R} y_r \leq U
&&
\label{BAND-BLP_activationBound}
\\
&x^{bs}_{ij} \in \{0,1\}
&& b \in B, s \in S, (i,j) \in A
\nonumber
\\
&y_r \in \{0,1\}
&& r \in R \; .
\nonumber
\end{align}

\noindent
We remark that the constraints (\ref{BAND-BLP_bioConservation}-\ref{BAND-BLP_sinkConservation}) could be simplified by dividing both sides of the inequalities by $d_{bs}$.

The objective function pursues the minimization of the total BAN energy consumption expressed as the sum of the energy consumed by each arc $(i,j)$, equal to the product of the data flow  and the energy $E_{ij}$ consumed on $(i,j)$ to transmit and receive 1 unit of data.
The constraints (\ref{BAND-BLP_bioConservation}-\ref{BAND-BLP_sinkConservation}) are flow conservation constraints essentially expressing the balance between ingoing and outgoing flows in any node of the graph. Note that we distinguish three flow balance cases, one for each type of device/vertex: 1) biosensors $b \in B$, which only transmit data, have a negative flow balance; 2) relays $r \in R$, which are transit vertices and thus retransmit all the received data, have a null flow balance; 3) sinks $s \in S$, which only receive data, have a positive flow balance. In each of these vertices, the flow balance must be considered for the data flow generated by each biosensor $b \in B$ for each sink $s \in S$.
The constraints \eqref{BAND-BLP_capacity} express the capacity of each relay. We note that each of these constraints has a right-hand-side whose value may vary: if $y_r = 1$, then the constraint activates and the right-hand-side is equal to $c_r$. Otherwise, the right-hand-side is equal to 0 and forces to zero also the left-hand-side, thus preventing data flows to be received by or transmitted to $r$.
The constraints \eqref{BAND-BLP_activationBound} express the limit $U > 0$ on the number of deployable relays

\smallskip
\noindent
\textbf{Protecting against data uncertainty by a robust model.}
\label{sec:ROB-BAND}
Until now, we have assumed that all the data involved in the BAND are exactly known when the problem is solved. However, this assumption does not hold in the real world: among its sensors, a BAN typically includes sensors that generate data according to an event-driven policy, thus leading to changeful and non-continuous data rates whose value is not known a priori \cite{ChEtAl10}.
A reduction in the expected data rate is not harmful to the designed BAN, since there is anyway sufficient transmission capacity. What can instead have very bad effects is an increase in the data rates, since the used relays may become not sufficient to handle the increased data volumes. In this case, biomedical data would be lost and this is a risk that cannot be taken at all in a BAN (as an example, one can think about the dramatic consequences that losing data produced by an early detection ischemia sensor could have on a patient).

The presence of data uncertainty in an optimization problem, namely the fact that a subset of the input data is not exactly known when the problem is solved, may result really tricky not just practically but also theoretically: as well-known from sensitivity analysis, even small deviations in the value of the input data may completely compromise the feasibility and optimality of produced solutions. Solutions supposed to be feasible may result infeasible and thus totally useless in practice, while solutions supposed to be optimal may result instead of very bad quality. For an exhaustive introduction to the consequences of the presence of data uncertainty in optimization, we refer the reader to \cite{BaEtAl14,BeElNe09,BuDA12a}.

For the design of BAN under data rate uncertainty, we adopt \emph{min-max robustness} (Min-Max) \cite{AiEtAl09,FuEtAl15}:
this type of robust optimization paradigm is especially appropriate in problems where it is crucial to guarantee very high level of protection against data uncertainty, since infeasibility due to data variations could have dramatic effects. This is the case of BANs, where data loss due to unexpected fluctuations in data rates may lead to the death of monitored patients.

In the context of the BAND, assuming the perspective of a highly risk-averse decision maker, who wants to guarantee a fully trustable monitoring of health conditions even under data uncertainty, looks appropriate.
Specifically, Min-Max can be adapted to the BAND problem by introducing a set of data generation scenarios $\Sigma$: each scenario $\sigma \in \Sigma$ specifies a vector $d^{\sigma} = (d^{\sigma}_{11}  \cdots d^{\sigma}_{bs} \cdots d^{\sigma}_{|B||S|})$ that states the bitrate between each biosensor-sink couple in $\sigma$. These scenarios can be included in the so-called \emph{robust counterpart} of BAND-ILP, a version of BAND-ILP that we denote by Rob-BAND-ILP and that produces robust solutions, i.e. solutions protected against data fluctuations (such formulation has been introduced in \cite{DANa15}, paper to which we refer the reader for a detailed discussion of how the model is derived following the principles of Min-Max):
\begin{align}
\min \hspace{0.3cm}
& E
&&
\mbox{(Rob-BAND-ILP)}
\label{Rob-BAND-BLP_objFunction}
\\
&
\sum_{b \in B} \sum_{s \in S} \sum_{(i,j) \in A}
\hspace{0.1cm} E_{ij} \hspace{0.1cm} d_{bs}^{\sigma} \hspace{0.1cm} x^{bs}_{ij}
\leq E
&&
\sigma \in \Sigma
\label{Rob-BAND-BLP_UBobjFunction}
\\
&
- \sum_{\substack{(b,j) \in \\ A_{B \rightarrow R} \cup A_{B \rightarrow S}}} x^{bs}_{bj} \hspace{0.1cm} = - 1
&&
b \in B, s \in S
\label{Rob-BAND-BLP_bioConservation}
\\
&
\sum_{\substack{(r,j) \in \\ A_{R \leftrightarrow R} \cup A_{R \rightarrow S}}} x^{bs}_{rj}
- \sum_{\substack{(j,s) \in \\ A_{B \rightarrow R} \cup A_{R \leftrightarrow R}}} x^{bs}_{jr}
 = 0
&&
b \in B, s \in S, r \in R
\label{Rob-BAND-BLP_relayConservation}
\\
&
\sum_{\substack{(j,s) \in \\ A_{B \rightarrow S} \cup A_{R \rightarrow S}}} x^{bs}_{js} \hspace{0.1cm} = \hspace{0.1cm} 1
&&
b \in B, s \in S
\label{Rob-BAND-BLP_sinkConservation}
\\
&
\sum_{\substack{(r,j) \in \\ A_{R \leftrightarrow R} \cup A_{R \rightarrow S}}} d_{bs}^{\sigma} \hspace{0.1cm} x^{bs}_{rj} \hspace{0.1cm} \leq \hspace{0.1cm} c_r \hspace{0.05cm} y_r
&&
r \in R, \sigma \in \Sigma
\label{Rob-BAND-BLP_capacity}
\\
&
\sum_{r \in R} y_r \leq U
&&
\label{Rob-BAND-BLP_activationBound}
\\
&x^{bs}_{ij} \in \{0,1\}
&& b \in B, s \in S, (i,j) \in A
\nonumber
\\
&y_r \in \{0,1\}
&& r \in R \; .
\nonumber
\end{align}

\noindent
where variable $E \geq 0$ is introduced to express an upper bound on the total energy consumed by routing decisions over all the scenarios in $\Sigma$. Additionally, this robust model includes one capacity constraint (\ref{Rob-BAND-BLP_capacity}) and one variable lower bound constraint \eqref{Rob-BAND-BLP_UBobjFunction} for each scenario $\sigma \in \Sigma$, as done in \cite{DANa15}.

\section{A fast heuristic for the Rob-BAND-ILP}
\label{sec:ANTS}

\noindent
The optimization problem Rob-BAND-ILP in principle can be solved by any Mixed Integer Programming (MIP) solver. However, the problem may prove (very) difficult to solve even for commercial solvers based on state-of-the-art branch-and-cut solution algorithms like IBM ILOG CPLEX \cite{CPLEX}: these solvers can have issues in fast finding good quality solutions and tend to present a slow convergence to an optimal solution. More specifically, the optimality gap, expressing how far the best solution found is from an optimal solution during the execution of the branch-and-cut, tends to be improved slowly. Such difficulties constitute an issue for practical applications of the BAND.

To tackle such unsatisfactory performance, we propose to adopt a heuristic that is based on the sequential execution of the following three phases:\\
1) a deterministic variable fixing phase that exploits the optimal solution of the (strengthened) linear relaxation of Rob-BAND-ILP and that produces a partial solution for the problem;\\
2) a probabilistic variable fixing phase, guided by the combination of information coming from the optimal solutions of two distinct linear relaxations of (Rob-)BAND-ILP and that provides a complete fixing of the variables;\\
3) a reparation/improvement phase based on executing an \emph{exact large variable neighborhood search}, which aims at substituting an infeasible fixing with a feasible fixing or improving a feasible fixing, produced during phase 2. The search is called \emph{exact} since it is expressed through the solution of an integer linear programming problem solved through an MIP solver.

We detail the features of each phase in the following subsections. Here, we just anticipate that the second phase represents the core of our algorithm and is based on an improvement of the algorithm ANTS (\emph{Approximate Nondeterministic Tree Search}) \cite{Ma99}, an ant colony-like algorithm. Specifically, the refinement of ANTS that we adopt is based on interpreting ant colony as a probabilistic variable fixing procedure, where the fixing is guided by optimal solutions to linear relaxations of the problem. Such interpretation has been first made in \cite{DAKrPu14,DAKrPu15} and \cite{DANa15}. We stress that such interpretation actually leads to an algorithm that in spirit and substance is deeply different from ant colony algorithms and is more ``well founded'' on precious polyhedral considerations that come from the linear relaxation of the problem and from its strengthened version obtained after the application of cuts at the root of the branch-and-bound node.

Ant Colony Optimization (ACO) is a metaheuristic inspired by the foraging behaviour of ants that was initially proposed for combinatorial optimization by Dorigo and colleagues in \cite{DoMaCo96} and later extended and improved in many works (e.g., \cite{GaMoWe12,Ma99} - we refer the reader to \cite{Bl05} for an overview of theory and applications of ACO).
A typical ACO algorithm has the following general structure:
\begin{description}
  \item \textbf{while} an arrest condition is \emph{not} satisfied
  \begin{description}
      \item - ant-based solution construction
      \item - pheromone trail update
  \end{description}
  \item - local search
\end{description}

The core of the algorithm is represented by a cycle where a number of feasible solutions are defined in a probabilistic and iterative way, taking into account the quality of solutions built in previous cycle iterations.
Each solution is iteratively built by an \emph{ant}: at each iteration, the ant is in a \emph{state} that corresponds to a \emph{partial solution} and can execute a so-called \emph{move}, fixing the value of an additional variable and thus further completing the solution. The move is established probabilistically, putting together an \emph{a-priori} and an \emph{a-posteriori} measure of variable fixing attractiveness. In the theory of ACO, the a-priori measure is called \emph{pheromone trail value} and is updated at the end of the construction phase on the basis of how good were the moves done. The construction cycle ends when reaching a stop condition, which commonly consists in a time limit. Then a local search is started to try improve the feasible solutions built, finding some local optimal solution.
\\
\textbf{Deterministic fixing phase.}
The first step of the first phase consists of solving the linear relaxation of problem Rob-BAND-ILP: the optimal solution of the linear relaxation, strengthened by the cuts added by CPLEX at the root node of the branch-and-bound tree, is then used to fix the value of decision variables of the problem. Specifically, the strategy is to fix to 1 the relay activation variables $y_r$ whose value in the optimal solution of the linear relaxation is sufficiently close to 1. The rationale at the basis of this strategy is that if the value of a variable is sufficiently close to 1, then there is a quite good indication that we should fix the decision variable to 1 in a good feasible solution.
Formally, if we denote by $y_r^{\text{TLR}}$ the value of variable $y_r$ in an optimal solution of the (strengthened) linear relaxation, the fixing rule is: if $y_r^{\text{TLR}}$ $\geq 1 - \epsilon$ then impose $y_r = 1$, where $\epsilon > 0$ is a parameter to choose. We focus on the relay activation variables as their fixing results particularly effective in reducing the difficulty of solving the complete problem Rob-BAND-ILP.
Once that this fixing has been operated, we obtain a smaller version of the original problem, denoted by Rob-BAND-ILP$^{FIX}$, where the fixed variables $y_r$ are no more part of the decision problem.
\\
\textbf{Probabilistic fixing phase.}
This phase is aimed at identifying the data routing paths within the BAN and consists of fixing the unsplittable flow variables $x_{ij}^{bs}$.
As first step, let us denote by $C$ the set of biosensor-sink pairs for which there exists at least one data scenario with positive bitrate, i.e. $C = \{(b,s) \in  B \times S: \exists \hspace{0.05cm} \sigma \in \Sigma \mbox{ with } d_{bs}^{\sigma} > 0\}$. We then refer to the concept of routing state.
\begin{definition}[Routing state - RS \cite{DANa15}]
  Consider a subset of biosensor-sink couples $\bar{C} \subseteq C$. We define \emph{routing state} a fixing of the unsplittable flow variables $x_{ij}^{bs}$ $\forall (i,j) \in E$ for each $(b,s) \in \bar{C}$ such that the fixing is feasible for the flow conservation constraints (\ref{Rob-BAND-BLP_bioConservation}-\ref{Rob-BAND-BLP_sinkConservation}).
\qed
\end{definition}

\noindent
A \emph{routing state} assigns one routing path to each pair $(b,s) \in \bar{C}$. It is said \emph{partial} when $\bar{C} \subset C$ (i.e., only a subset of data flows is routed), whereas it is said \emph{complete} when $\bar{C} = C$ (i.e., all data flows are routed).

We build a complete routing state by assigning paths to biosensor-sink pairs  according to the following procedure, which we call {\sc set-paths}. The pairs are considered according to an a-priori defined order, as done in \cite{DANa15}: we sort pairs $(b,s) \in C$ for decreasing value of the highest bitrate $d_{bs}^{\sigma}$ over all the data scenarios $\sigma \in \Sigma$.
Following the pair order, for each couple $c = (b,s) \in C$, we assign the entire data flow to a path $p$ connecting $b$ to $s$. The routing path for the pair $c$ is chosen from a set of candidates $P_c$, defined as follows: 1) we solve the linear relaxation of Rob-BAND-ILP$^{FIX}$, which includes the deterministic fixing of the first phase and where we have additionally fixed the value of variables of pairs $c \in \bar{C}$ for which a path has been assigned in previous executions of the external loop; 2) using the optimal solution $(x^{LR},y^{LR})$ of the linear relaxation, we define a graph $H^{c}(V,A^{mod})$ from $G(V,A)$: the set of vertices does not change, while in $A^{mod}$ we keep only those arcs $(i,j) \in A$ with a positive flow, i.e. such that $x^{LR \hspace{0.05cm} c=(b,s)}_{ij} > 0$. Furthermore, for each arc $(i,j) \in A^{mod}$ we define a  weight $w_{ij} = x^{LR \hspace{0.05cm} c=(b,s)}_{ij}$. We derive $L$ candidate paths for the pair $c = (b,s)$ on graph $H^{c}(V,A^{mod})$ by iteratively modifying $H^{c}(V,A^{mod})$: in an internal loop, at each iteration we find the shortest path $p$ considering the weights $w_{ij}$ in $H^{c}(V,A^{mod})$, then we add $p$ to the set $P_c$ and we delete the arc of $p$ with lowest weight from $H^{c}(V,A^{mod})$. This is a straightforward procedure that, however, can be fast implemented and that we have observed among professionals in real-world telecommunication applications. The rationale behind the exclusion of the arc with lowest weight is that, if the fractional value in the range [0,1] of a binary variable is seen as the probability of fixing to 1 the variable in a good solution, then smaller values should lead to fixing to 1 of lower quality (we refer the reader to the book \cite{MoRa95} about randomized rounding algorithms for a good discussion on looking at fractional binary solutions as measures of probability). After having established the set of candidate paths $P_c$ for $c$, we compute the probability of choosing each path $p \in P_c$ to route the entire flow of couple $c$ using the formula:
\begin{equation}
\label{ProbANTS}
\text{PROB}_{p} = \frac{\alpha \hspace{0.1cm} \tau_{p}  + (1-\alpha) \hspace{0.1cm} \eta_{p} }
                {\sum_{\pi \in P_c} \alpha \hspace{0.1cm} \tau_{\pi}  + (1-\alpha) \hspace{0.1cm} \eta_{\pi} }
\; \;,
\end{equation}
where both $\tau_{p}$ and $\eta_{p}$ are obtained as the sum of the current values of the a-priori and the a-posteriori measures $\tau_{ij}^{c}$, $\eta_{ij}^{c}$ for the edges in path $p$ for pair $c$. In particular, the a-priori measures $\tau_{ij}^{c}$ are initialized with the values that flow variables assume in an optimal solution to the (strengthened) linear relaxation of Rob-BAND-ILP$^{FIX}$ and are updated at the end of each construction phase (see below for more details). Instead, the a-posteriori measures $\eta_{ij}^{c}$ are set equal to  the values that flow variables assume in an optimal solution to the linear relaxation of Rob-BAND-ILP$^{FIX}$ plus the additional fixing that have been operated while building a complete routing state.
After having probabilistically chosen a path $p^{*} \in P_c$ through formula \eqref{ProbANTS}, we derive a fixing of the flow variables $x^{bs}$, where $x^{bs}_{ij} = 1$ if $(i,j)$ belongs to $p^{*}$ and $x^{bs}_{ij} = 0$ otherwise. Finally, we add the couple $c$ to the set of processed couples $\bar{C}$ for which the routing has been established.

After having executed  the external loop $|C|$ times, following the ordering of the pairs, we obtain a complete routing state.
However, since the procedure adopted to define a routing state does not take into account the capacity of relays, we may actually produce routing solutions that are infeasible: this can occur, for example, if many routing paths use the same relay and the sum of the data exceeds the relay capacity. Due to this possibility, we include in the algorithm a check-and-repair phase: this phase first verifies the feasibility of the routing state and, in case of infeasibility, tries to repair the solution. The reparation is attempted through the same ILP heuristic that we adopt to find better solutions (see the next subsection for a description of the ILP heuristic).

The feasibility of a complete routing state for the complete problem Rob-BAND-ILP can be fast and easily operated: we deploy all relays appearing in paths used in the routing state (i.e., we fix to 1 the corresponding relay deployment variables $y_r$ and to 0 all the other variables $y_r$) and we verify the presence of relay-capacity constraints \eqref{Rob-BAND-BLP_capacity} violated for some data scenario in $\Sigma$. Additionally, we must check if the number of activated relays is higher than the limit expressed by constraint \eqref{Rob-BAND-BLP_activationBound}. If all constraints are satisfied, then we have built a feasible solution for Rob-BAND-ILP: the complete routing state specifies the values of the flow variables $x$ and
these allows us to also derive a feasible activation of the relay activation variables $y$. In contrast, any violation of a capacity or activation constraint immediately certifies that the built routing state is infeasible and we must therefore repair it.

We present the pseudocode of the heuristic in Algorithm \ref{ALG_RobuBAND}. There, the energy value of a solution $(\bar{x},\bar{y})$ is denoted by $E(\bar{x},\bar{y})$. Additionally, we denote by  $(x^{*},y^{*})$ the best solution found during the entire execution of the algorithm. The heuristic includes two main loops: the external loop is executed until a time limit is reached, whereas the internal loop  provides for building $m$ feasible solutions according to the routing state construction that we have explained above. Specifically, the first task of the algorithm is to solve the (strengthened) linear relaxation of Rob-BAND-ILP that is used to execute the first deterministic fixing phase, leading to problem Rob-BAND-ILP$^{FIX}$. Then the (strengthened) linear relaxation of Rob-BAND-ILP$^{FIX}$ is solved and its optimal solution is used to initialize the a-priori measure of attractiveness $\tau_{ij}^{c}(0)$.
In each execution of the internal loop, the first task is to define a complete routing state as previously detailed. The complete routing state provides a complete valorization of the variables $\bar{x}$ and is used as basis to derive a relay installation $\bar{y}$. This leads to an integral solution $(\bar{x},\bar{y})$ whose feasibility is not guaranteed and must thus be checked and eventually repaired through the ILP heuristic. If the solution $(\bar{x},\bar{y})$ found is feasible and is better than the current best solution $(x^{B},y^{B})$, $(x^{B},y^{B})$ is updated and the internal loop continues.
At the end of the internal loop,
the a-priori measures $\tau$ are updated, evaluating how good the fixing resulted in the obtained solutions.
The update formula uses the \emph{optimality gap} (\emph{OGap}) corresponding with a feasible solution of value $v$ and a lower bound $L$ for the optimal value $v^{*}$ of the problem (since we consider a minimization problem, we have $L \leq v^{*} \leq v$ and $OGap(v,L)$ $= (v - L)/v$):
\begin{equation}
\small
\label{updateFormula}
\tau_{ij}^{c}(h) = \tau_{ij}^{c}(h-1)
+ \sum_{f=1}^{F} \Delta \tau_{ij}^{c}
\mbox{ with }
\Delta \tau_{ij}^{c \; f} =
\tau_{ij}^{c}(0)
\cdot \left(
\frac{OGap(\bar{v},L) - OGap(v_{f},L)}{OGap(\bar{v},L)}
\right)
\end{equation}

\noindent
where $\tau_{ij}^{c}(h)$ is the a-priori attractiveness of fixing variable $x_{ij}^{c = (b,s)}$ at fixing iteration $h$, $L$ is the lower bound (in our case, the strengthened linear relaxation of Rob-BAND-ILP),
$v_{f}$ is the value of the $f$-th feasible solution built in the last construction cycle and $\bar{v}$ is the (moving) average of the values of the $F$ solutions produced in the previous construction phase. $\Delta \tau_{ij}^{c \; f}$ is the penalization/reward factor for a fixing and depends upon the initialization value $\tau_{ij}^{c}$ of $\tau$, combined with the relative variation in the optimality gap that $v_{\sigma}$ implies w.r.t. $\bar{v}$.
Since in \eqref{updateFormula} we use a relative gap difference, we are able to encourage or discourage fixings made in the last produced solution through a comparison with the average quality of the last $F$ solutions produced.
Once the time limit is reached, we execute the ILP heuristic for improving the best solution found and at the end of the execution we return $(x^{*},y^{*})$.
\begin{algorithm}
\caption{Heuristic for Min-Max BAND}
\label{ALG_RobuBAND}
\begin{algorithmic}[1]
    \State compute the strengthened linear relaxation of Rob-BAND-ILP
    \State execute the deterministic fixing phase of variables $y_r$ using a fixing threshold $\epsilon > 0$ and define Rob-BAND-ILP$^{FIX}$
    \State compute the strengthened linear relaxation of Rob-BAND-ILP$^{FIX}$ and initialize the values $\tau_{ij}^{c}(0)$ through it
    \State let $(x^{*},y^{*})$ denote the best solution found by the algorithm
    \While{a global time limit is not reached}
        \State let $(x^{B},y^{B})$ denote the best solution found in the inner loop
        \For{$k := 1$ to $m$}
            \State build a complete routing state $\bar{x}$ following the procedure {\sc set-paths}
            \State derive a relay installation $\bar{y}$ using $\bar{x}$
            \If {$(\bar{x},\bar{y})$ is not feasible for Rob-BAND-ILP}
                \State run mod-RINS for repairing $(\bar{x},\bar{y})$
            \EndIf
            \If {$(\bar{x},\bar{y})$ is feasible and $E(\bar{x},\bar{y}) < E(x^{B},y^{B})$}
                \State update the best solution found $(x^{B},y^{B}) := (\bar{x},\bar{y})$
            \EndIf
        \EndFor
        \State update $\tau_{ij}^{c}(t)$ according to (\ref{updateFormula})
        %
        \If {$E(x^{B},y^{B}) < E(x^{*},y^{*})$}
            \State update the best solution found $(x^{*},y^{*}) := (x^{B},y^{B})$
        \EndIf
    \EndWhile
    \State run mod-RINS($x^{*},y^{*}$) for improving $(x^{*},y^{*})$
    \State return $(x^{*},y^{*})$
\end{algorithmic}
\end{algorithm}

\noindent
\textbf{Reparation/improvement by an ILP heuristic.}
To either repair an infeasible fixing of the decision variables or to improve a  feasible solution, we rely on an ILP heuristic \emph{exactly} executing a \emph{very large neighborhood search}, i.e. the search is formulated as an integer linear programming problem solved through a state-of-the-art MIP solver (see also, e.g., \cite{DACa16}). Specifically, we rely on a modified \emph{Relaxation Induced Neighborhood Search} (RINS) (see \cite{DaRoLP05} for an exhaustive description of this search method).
Let $(\bar{x},\bar{y})$ be a solution found for Rob-BAND-ILP, $(x^{\small TLR},y^{\small TLR})$ be an optimal solution of the linear relaxation of Rob-BAND-ILP$^{FIX}$ strengthened by the cuts found by CPLEX in the root node of the branch-and-bound tree.
Moreover, denote by $(\bar{x},\bar{y})_j, (x^{\small TLR},y^{\small TLR})_j$ the $j$-th component of the vectors.

The modified RINS \emph{(mod-RINS)} that we adopt solves a subproblem of Rob-BAND-ILP where we fix the variables whose value in $(\bar{x},\bar{y})$ and $(x^{\small TLR},y^{\small TLR})$ differs of at most $\rho > 0$ according to the following rules:\\
\; $(\bar{x},\bar{y})_j = 0$
          $\hspace{0.1cm} \wedge \hspace{0.1cm}$
          $(x^{\small TLR},y^{\small TLR}) \leq \rho$
          $\hspace{0.2cm} \Longrightarrow \hspace{0.2cm}$
          $(x,y)_j = 0$\\
\; $(\bar{x},\bar{y})_j = 1$
          $\hspace{0.1cm} \wedge  \hspace{0.1cm}$
          $(x^{\small TLR},y^{\small TLR}) \geq 1 - \rho$
          $\hspace{0.2cm} \Longrightarrow \hspace{0.2cm}$
          $(x,y)_j = 1$ \; .
\\
The resulting problem is then passed to CPLEX, which attempts at solving it within a time limit. The rationale
is that CPLEX, though not being able to fast finding good quality solutions to the complete problem, is instead able to fast finding good solutions to subproblems obtained by smartly fixing variables.

\section{Experimental results}
\label{sec:computations}

We evaluated the performance of the new heuristic on the same set of 30 instances considered in \cite{DANa15}. We refer the reader to that paper for a detailed description of the instances; here we just remind the major topology features of the corresponding graph: all instances consider a BAN including 16 biosensors (i.e., $|B| = 16$) and 2 sinks (i.e., $|S| = 2$). Moreover, 400 potential sites over the human body (excluding head, hands and feet) are considered for the deployment of relays (i.e., $|R| = 400$), chosen randomly over the human body.

We performed all the experiments on a 2.70 GHz machine with 8 GB, using a C/C++ code interfaced with IBM ILOG CPLEX 12.5 through Concert Technology and running with a time limit of 2400 second. The results of the computational tests are presented in Table \ref{table:results}, where \emph{ID} is the identifier of the instance and where we show the performance of all the considered algorithms in terms of the optimality gap associated with the best solutions found within the time limit.
Specifically, we show the optimality gaps of: 1) CPLEX ($\text{GapILP}\%$) applied directly to solve Rob-BAND-ILP; 2) the heuristic presented in \cite{DANa15} ($\text{GapRB}\%$) - we denote this heuristic by $\emph{RB}$; 3) our new heuristic ($\text{GapHEU}\%$) - we denote our heuristic by $\emph{HEU}$.
Finally, $\Delta Gap\%$ is the percentage increase of the optimality gap of CPLEX w.r.t. that of the heuristics. For both RB and HEU, the optimality gap is derived comparing the value of the linear relaxation of Rob-BAND-ILP computed by CPLEX with the value of the best feasible solution found by the heuristic within the time limit. For both heuristics, the number of candidate paths for each biosensor-sink routing path assignment equals 5, the combination factor of the a-priori and a-posteriori measures $\alpha$ is set equal to $0.5$ and the width $F$ of the moving average is $4$. To solve the linear relaxation of Rob-BAND-ILP and of BAND-ILP we used CPLEX.
The threshold $\epsilon$ for the deterministic fixing threshold is set equal to $10^{-1}.$
The repair/improvement heuristic \emph{mod-RINS} uses a threshold  $\rho = 10^{-1}$ and runs with a time limit of 10 minutes for finding improvements and of 1 minute when used for solution reparation. The external cycle of $\text{HEU}$ ran with a time limit of 30 minutes, which matches the time limit of CPLEX when added up to the time reserved for mod-RINS.

\begin{table}
\caption{Experimental results}
\label{table:results}
\begin{center}
\begin{tabular}{c | c c c| c c}
\hline
ID  & $\text{ GapILP}\% $ & $\text{ GapRB}\% $ 	& $\text{ GapHEU}\% $ \;	 &
\; $\Delta \text{Gap$\%$ HEU-ILP}$ \;	&  $\Delta \text{Gap$\%$ HEU-RB}$
\\ [2pt]
\hline
I1 	&  55.22	&  	41.63	& 	36.85 &  33.2 &  11.5\\
I2 	&  67.17 	&  	57.10	& 	49.28	&  26.6 &  13.7\\
I3 	&  40.26	&  	35.83	& 	32.21	&  19.9 &  10.1\\
I4 	&  45.20	&  	42.16	& 	26.10	&  25.1 &  10.1\\
I5 	&  68.60	&  	54.54	& 	51.65	&  24.7 &  5.3\\
I6 	&  60.45	&  	38.33	& 	35.34	&  41.5 &  7.8\\
I7 	&  45.65	&  	34.52	& 	33.07	&  27.5 &  4.2\\
I8 	&  64.09	&  	48.78	& 	44.68	&  30.2 &  8.4\\
I9 	&  60.66	&  	47.77	& 	42.13	&  30.5 &  11.8\\
I10 	&  34.08	&  	28.10	& 	23.07	&  32.3 &  17.9\\
I11 	&  61.42	&  	48.50	& 	32.20	&  46.2 &  14.5\\
I12 	&  60.97	&  	46.65	& 	43.56	&  31.8 &  14.2\\
I13 	&  59.96	&  	37.66	& 	27.27	&  21.2 &  5.8\\
I14 	&  63.92	&  	50.77	& 	28.73	&  25.0 &  9.9\\
I15 	&  34.61	&  	28.95	& 	32.47	&  10.9 &  21.0\\
I16 	&  38.33	&  	31.89	& 	40.21	&  34.5 &  17.1\\
I17 	&  36.45	&  	41.10	& 	41.23	&  32.3 &  11.6\\
I18 	&  53.04	&  	41.61	& 	32.17	&  39.3 &  22.7\\
I19 	&  36.81	&  	31.97	& 	28.36	&  22.9 &  11.3\\
I20 	&  36.08	&  	30.87	& 	26.70	&  25.9 &  13.5\\
I21 	&  34.89	&  	29.03	& 	35.96	&  20.4 &  14.7\\
I22 	&  56.89	&  	42.63	& 	37.00	&  34.9 &  13.2\\
I23 	&  52.83	&  	43.58	& 	39.92	&  24.4 &  8.4\\
I24 	&  47.13	&  	50.34	& 	40.93	&  13.1 &  18.7\\
I25 	&  41.50	&  	35.06	& 	37.34	&  10.0 &  6.5\\
I26 	&  67.46	&  	38.43	& 	35.47	&  47.4 &  7.7\\
I27 	&  34.18	&  	29.52	& 	24.38	&  28.6 &  17.4\\
I28 	&  64.75	&  	52.11	& 	40.48	&  37.4 &  22.3\\
I29 	&  70.31	&  	41.25	& 	34.61	&  50.7 &  16.1\\
I30 	&  59.95	&  	48.35	& 	42.50	&  29.1 &  12.1\\
\hline
\end{tabular}
\end{center}
\end{table}
The optimality gaps $\text{GapILP}\%$ indicate that the instances proved challenging to solve even for a state-of-the-art solver like CPLEX, which produces solutions that are distant from the optimum at the time limit. In contrast, HEU provides solutions associated with a great reduction in the optimality gap that is on average equal to about 29\%  and can be sensibly over 40\%, as in the case of instances I11, I26 and I29. HEU is also able to grant a significative reduction with respect to the benchmark heuristic RB, which already grants a high reduction in the optimality gap with respect to CPLEX: the average reduction in gap is about 12\% and can be over 20\% as for I18 and I28.
We think that the better performance of HEU with respect to BR is due to the inclusion of an additional fixing phase that involve the relay activation variables, which are excluded from the linear relaxation-guided fixing of BR.

As future work, we plan to further reduce the optimality gap by considering other integration of heuristics (in particular, genetic and sequential heuristics like in \cite{DA11,DeDAKa15}) and cutting plane methods identifying conflicts between variables, as in \cite{BlDAKa13,DAMaSa11}.
Also, we intend to evaluate biobjective versions of the problem, considering the trade-off between relay deployment cost and energy consumption and adopting an algorithm similar to \cite{ZaEtAl13}. Finally, we plan to investigate the adoption of another robustness paradigm, namely Multiband Robustness \cite{BuDA12a}.

\bibliographystyle{splncs}
\bibliography{DAndreagiovanni_EvoStar2017}

\end{document}